\documentclass[10pt,reqno]{amsart}

\usepackage{imports}

\title{Threshold existence for random NAE-$k$-SAT and random hypergraph $2$-colourability}

\author{Pavel~Zakharov}

\address{Pavel Zakharov, {\tt pavel.zakharov@uni-bielefeld.de}, Bielefeld University, Faculty of Mathematics, Universit\"atsstra\ss e~25, 33615 Bielefeld, Germany.}

\thanks{
	Pavel Zakharov is supported by DFG CO 646/6.
}

\begin{document}

\begin{abstract}
	The existence of limiting threshold densities for random hypergraph
	$2$-colourability and random NAE-$k$-SAT was conjectured in earlier
	work, including Achlioptas and Moore [RANDOM 2002] and Coja-Oghlan
	and Panagiotou [STOC 2012].
	We prove these conjectures for all $k\geq 3$.
	We also show that $\alpha^{\mathrm{2-col}}_k \leq
	\alpha^{\mathrm{NAE}}_k$ for each $k\geq 3$.
\end{abstract}
\maketitle

\section{Introduction}

For the study of random constraint satisfaction problems (CSPs), one of the starting points was the question of \Erdos\ and \Renyi\ about the chromatic number of a random graph~\cite{ER}.
In the language of random CSPs, fixing the number of colours leads to the question of the $q$-colourability threshold.
Estimates of the chromatic number are known in many regimes \cite{AN, CIKZ, Heckel}, but the existence of a limiting colourability threshold remains open for general $q$.
One of the explanations regarding why the problem is so difficult comes from statistical mechanics.
A heuristic called the cavity method predicts that before the satisfiability threshold, the solution set of a random CSP can undergo several other phase transitions, such as clustering and condensation~\cite{KMRTZ}.
Among other things, these transitions help explain the limitations of the second moment method~\cite{CP, CZ}, one of the main tools for estimating thresholds in random structures.

A related question is the satisfiability conjecture for random $k$-SAT~\cite{CR}.
Here a formula consists of independently sampled clauses, each being a disjunction of $k$ distinct uniformly chosen variables with independent uniform signs, and an assignment satisfies the formula if it makes at least one literal in every clause true.
The conjecture asserts that for each fixed $k\geq3$ there is a constant clause density $\alpha_k^{\mathrm{SAT}}$ below which such a formula is satisfiable with high probability (\whp) and above which it is unsatisfiable \whp.

Friedgut's sharp threshold theorem~\cite{Fri} establishes the existence of a threshold sequence, but does not establish its convergence.
Writing $\Phi^{\mathrm{SAT}}_{n,m}$ for the random $k$-SAT formula described above, for every fixed $k\geq3$ there is a sequence $\alpha_k^{\mathrm{SAT}}(n)$ such that for every fixed $\eps>0$ and every sequence $\alpha=\alpha(n)\geq0$,
\begin{align*}
	\Pr\bc{\Phi^{\mathrm{SAT}}_{n,\floor{\alpha n}}\text{ is satisfiable}}
	\rightarrow \begin{cases}
		1, & \text{ if $0 \leq \alpha < \alpha_k^{\mathrm{SAT}}(n) \cdot (1- \varepsilon)$};
		\\
		0, & \text{ if $\alpha > \alpha_k^{\mathrm{SAT}}(n) \cdot (1+ \varepsilon)$},
	\end{cases}
\end{align*}
as $n$ goes to infinity, provided the corresponding inequality holds for all sufficiently large $n$.
Similar sharp threshold results also hold for graph colourability~\cite{AF} and for hypergraph $2$-colourability and NAE-$k$-SAT~\cite{FriH}.
However, sharpness alone does not imply that the threshold sequences converge as $n\ra\infty$.

For random $2$-SAT the threshold density is known to be $\alpha^{\mathrm{SAT}}_2 = 1$~\cite{CR}.
For larger clause sizes, Ding, Sly and Sun~\cite{DSS} proved the satisfiability conjecture for all sufficiently large $k$ and identified the threshold with the value predicted by the cavity method.
Recently, Carenini established an $O(n/\log n)$ scaling window~\cite{Car2}, and then improved it to $O(n^{1/2+1/k})$~\cite{Car}.
These bounds count the number of random clauses one needs to add to decrease the satisfiability probability from $1-\eta$ to $\eta$ for any fixed $0<\eta<1/2$.
Combined with the Abbe--Montanari interpolation criterion~\cite{AbMo}, the latter estimate proves the existence of a limiting threshold for every fixed $k\geq3$ without identifying its value.
An alternative proof later appeared in a manuscript by OpenAI~\cite{AI}.
It combines polynomial concentration for the last satisfiable index with an interpolation method to prove convergence directly, without using the Abbe--Montanari criterion.

For random hypergraph $2$-colourability and random NAE-$k$-SAT, the existence of limiting threshold densities was conjectured in~\cite{AM,CP}.
Let $\alpha_k^{\mathrm{2-col}}(n)$ and $\alpha_k^{\mathrm{NAE}}(n)$ denote the corresponding sharp threshold sequences.
For all sufficiently large $k$ and $n$, the bound of Coja-Oghlan and Zdeborov\'a~\cite{CZ} and the first moment method give
\begin{align*}
	2^{k-1}\ln2-\ln2-o_k(1)
	\leq\alpha_k^{\mathrm{2-col}}(n)
	\leq2^{k-1}\ln2-{\ln2\over2}+o_k(1).
\end{align*}
For NAE-$k$-SAT, Coja-Oghlan and Panagiotou~\cite{CP} obtained the sharper estimate
\begin{align}\label{eq_nae_aco}
	\alpha_k^{\mathrm{NAE}}(n)=2^{k-1}\ln2-{\ln2\over2}-{1\over4}+o_k(1).
\end{align}
Here the errors tend to zero as $k\ra\infty$ uniformly for sufficiently large $n$, but these estimates do not establish convergence in $n$ for any fixed $k$.
Coja-Oghlan and Panagiotou also observe that the results of~\cite{CZ} carry over to NAE-$k$-SAT, and it is natural to expect the same asymptotic threshold for hypergraph $2$-colourability.
More recently, Carenini~\cite{Car2} proved $O(n/\log n)$ scaling windows for both models, again leaving the convergence question open.

In this paper we prove that the limiting densities exist in both models for every fixed $k\geq3$.
For $n\geq k$, let $H_{n,m}$ be the random $k$-uniform hypergraph obtained by sampling $m$ edges independently and uniformly with replacement from the $k$-subsets of its vertex set.
A NAE-clause is satisfied if its literal values are not all equal.
We represent $H_{n,m}$ by a NAE-formula $\Phi^+_{n,m}$ with all signs positive, so the formula is satisfiable precisely when the hypergraph is $2$-colourable. Let $\Phi^\pm_{n,m}$ be the corresponding random NAE-$k$-SAT formula with independent uniform signs.
Vertices correspond to variables, edges to clauses and colourings to assignments. In particular, variables within each clause are distinct, while clauses may repeat.
We write $\SAT$ for the set of satisfiable NAE-formulas, and in statements concerning both models, $\alpha_k^\star$ denotes $\alpha_k^{\mathrm{2-col}}$ for $\star=+$ and $\alpha_k^{\mathrm{NAE}}$ for $\star=\pm$.

\begin{theorem}\label{thm_main}
	For every fixed $k\geq3$ there exist constants $\alpha_k^{\mathrm{2-col}}$ and $\alpha_k^{\mathrm{NAE}}$ such that for each $\star\in\{+,\pm\}$ and $m=\floor{cn}$,
	\begin{align*}
		\lim_{n\ra\infty}\Pr\bc{\Phi^\star_{n,m}\in\SAT}
		=\begin{cases}
			1, & 0\leq c<\alpha_k^\star;\\
			0, & c>\alpha_k^\star,
		\end{cases}
	\end{align*}
\end{theorem}
\noindent
Our second result compares the two thresholds for every clause size.
\begin{theorem}\label{thm_comparison}
	In the setup of \Thm~\ref{thm_main}, $\alpha_k^{\mathrm{2-col}} \leq \alpha_k^{\mathrm{NAE}}$.
\end{theorem}
Combining \Thm~\ref{thm_comparison} with \eqref{eq_nae_aco} gives the following upper bound for the random hypergraph $2$-colourability threshold.
\begin{corollary}
	\begin{align*}
		\alpha_k^{\mathrm{2-col}} \leq 2^{k-1} \ln 2 - {\ln 2 \over 2}  - {1 \over 4} + o_k(1).
	\end{align*}
\end{corollary}
\noindent
The proof of \Thm~\ref{thm_main} follows the approach developed for random $k$-SAT in~\cite{AI}.
For $\star\in\{+,\pm\}$, add corresponding clauses one by one, and let
\begin{align*}
	\xi^\star_n = \inf\cbc{m\geq1:~\Phi^\star_{n,m}\notin\SAT},\qquad
	\tau^\star_n=\min\cbc{\xi^\star_n-1,2^k n}.
\end{align*}
Thus $\tau_n^\star$ is the last satisfiable index, capped at $2^k n$.

We first show that for $k \geq 3$, $\tau_n^\star$ is within $n^{19/20}$ of its expectation with probability polynomially close to 1. A Poisson interpolation then compares satisfiability probabilities of the random formula on disjoint sets of variables with the random formula on their union, which is similar to the approach of Bayati, Gamarnik and Tetali~\cite{BGT}.
Together, these estimates imply that $\Erw\tau_n^\star/n$ converges.

The main novelty concerns the symmetry of the considered CSPs.
For NAE-satisfiability, flipping all variable values preserves a solution. Thus no variable has the same value in all solutions, and the argument based on frozen variables does not apply directly. Instead, we group variables according to whether their relative values are fixed across all solutions. These equivalence classes are used to determine exactly the probability that an additional clause destroys every solution, which yields both the clause replacement estimate needed for concentration, and the convexity inequality needed for the interpolation.
Interpolating between models with positive and randomly signed clauses also gives \Thm~\ref{thm_comparison}.

\subsection*{Declaration of AI use} ChatGPT was used to assist with proofreading, presentation and checking proofs.
It also suggested using a similar interpolation argument to establish \Thm~\ref{thm_comparison}.

\section{Preliminaries}
Throughout, the clauses are sampled independently with replacement. Within each proper clause, variables are distinct and uniformly sampled, and assignments take values in $\{0,1\}$. A clause is satisfied if its literal values are not all equal.

\subsection{Equivalence classes and clause addition}\label{sec_classes}
Let $\Phi$ be any satisfiable NAE-formula on a variable set $V$ of size $N$, with arbitrary signs and clause lengths. Write $\cS(\Phi)$ for its set of satisfying assignments. Define
\[
	x\sim y\quad\Longleftrightarrow\quad
	\sigma(x)\oplus\sigma(y)\text{ is constant over }\sigma\in\cS(\Phi).
\]
Transitivity follows from $\sigma(x)\oplus\sigma(z)=(\sigma(x)\oplus\sigma(y))\oplus(\sigma(y)\oplus\sigma(z))$, and reflexivity and symmetry are immediate. Denote its equivalence classes by $\cO_1,\ldots,\cO_J$.
Fix $\sigma_*\in\cS(\Phi)$ and put
\[
	\cO_{j,a}=\{y\in\cO_j:\sigma_*(y)=a\},\qquad a\in\{0,1\}.
\]
Within each part of one class the variables are equal in every solution, and variables in opposite parts of one class are opposite in every solution. Choosing another $\sigma_*$ can only interchange the two parts of each class.

For $2\leq r\leq N$, let $q_r^+(\Phi)$ and $q_r^\pm(\Phi)$ be the probabilities that adding a uniformly random $r$-clause on distinct variables makes $\Phi$ unsatisfiable, with all signs being positive or independent and uniformly random respectively. Write $\widehat q_r^+(\Phi)$ and $\widehat q_r^\pm(\Phi)$ when the variables are instead sampled independently with replacement. For an unsatisfiable $\Phi$ we set all four quantities to zero.

\begin{lemma}\label{lem_clause_classes}
	For satisfiable $\Phi$ and $2\leq r\leq N$,
	\begin{align}\label{eq_classes_distinct}
		q_r^\pm(\Phi)&=2^{1-r}\sum_j\frac{(|\cO_j|)_r}{(N)_r},&
		q_r^+(\Phi)&=\sum_j\sum_{a=0}^1\frac{(|\cO_{j,a}|)_r}{(N)_r}.
	\end{align}
	For sampling with replacement and $r\geq2$,
	\begin{align}\label{eq_classes_replacement}
		\widehat q_r^\pm(\Phi)&=2^{1-r}\sum_j\left(\frac{|\cO_j|}{N}\right)^r,&
		\widehat q_r^+(\Phi)&=\sum_j\sum_{a=0}^1\left(\frac{|\cO_{j,a}|}{N}\right)^r.
	\end{align}
	Here $(u)_r=u(u-1)\cdots(u-r+1)$.
\end{lemma}
\begin{proof}
	A new clause destroys every solution precisely when its literal values are equal in every solution. For random signs this requires all its variables to lie in one class $\cO_j$, and conditional on such a choice of variables, exactly two of the $2^r$ sign configurations have this property. This remains true with repeated variables, since signs are sampled independently for each literal. For positive signs the variables must all lie in one part $\cO_{j,a}$ of one class.
\end{proof}

The next lemma compares the probabilities of destroying all solutions by clauses of consecutive lengths.
\begin{lemma}\label{lem_consecutive}
	Let $2\leq r<N$ and let $\Phi$ be a satisfiable NAE-formula on $N$ variables, with arbitrary signs and clause lengths. For $\star\in\{+,\pm\}$ put $q=q_r^\star(\Phi)$. Then
	\begin{align}
		q_{r+1}^\star(\Phi)
		\geq q^{r/(r-1)}-\frac{r-1}{N-r}q.\nonumber
	\end{align}
\end{lemma}
\begin{proof}
	For $\star=\pm$, let $S_i=|\cO_i|$, with $i\in[J]$. For $\star=+$, let $S_i=|\cO_{j,a}|$, where $i=(j,a)$ ranges over the pairs in $[J]\times\{0,1\}$. Put $a_\ell^\pm=2^{1-\ell}$ and $a_\ell^+=1$. By \Lem~\ref{lem_clause_classes},
	\begin{align*}
		q_\ell^\star(\Phi)=a_\ell^\star\sum_i{(S_i)_\ell\over(N)_\ell}.
	\end{align*}
	Set $y_i=(S_i-1)_{r-1}/(N-1)_{r-1}$ for $S_i\geq r$, and $y_i=0$ for $S_i<r$. Then
	\begin{align}\label{eq_lem_consecutive_1}
		\sum_i{S_i\over N}=1,\qquad
		\sum_i{S_i\over N}y_i = {q \over a^\star_r},\qquad
		q_{r+1}^\star(\Phi)=a_{r+1}^\star\sum_i{S_i\over N}y_i{S_i-r\over N-r}.
	\end{align}
	For $S_i\geq r$ we have $y_i\leq\bc{(S_i-1)/(N-1)}^{r-1}$, and for $S_i<r$ we have $y_i=0$. Thus, in either case,
	\begin{align}\label{eq_lem_consecutive_2}
		y_i(S_i-r)\geq(N-1)y_i^{r/(r-1)}-(r-1)y_i.
	\end{align}
	Applying convexity of $x \mapsto x^{r/(r-1)}$ and using \eqref{eq_lem_consecutive_1}--\eqref{eq_lem_consecutive_2} yields
	\begin{align*}
		q_{r+1}^\star(\Phi)
		&\geq a_{r+1}^\star{N-1\over N-r}\bc{q\over a_r^\star}^{r/(r-1)}
		-{a_{r+1}^\star\over a_r^\star}{r-1\over N-r}q\\
		&\geq q^{r/(r-1)}-{r-1\over N-r}q,
	\end{align*}
	where the last step uses $a_{r+1}^\star / (a_r^\star)^{r/(r-1)}=1$ and $a_{r+1}^\star/a_r^\star\leq1$.
\end{proof}

\subsection{Forced variables}
We call a variable $x$ forced if $\Phi$ is satisfiable and flipping $x$ in any satisfying assignment makes that assignment unsatisfying. Let $b(\Phi)$ be the number of forced variables, with $b(\Phi)=0$ for unsatisfiable $\Phi$.
Every variable in a class of size at least two is forced. Consequently, for $2\leq r\leq N$ and $\star\in\{+,\pm\}$ \Lem~\ref{lem_clause_classes} implies
\begin{align}\label{eq_forced_bound}
	q_r^\star(\Phi)
	\leq\frac{(b(\Phi))_r}{(N)_r}\leq\left(\frac{b(\Phi)}N\right)^r.
\end{align}
\noindent
Finally, the following bounds keep the expected last satisfiable index away from zero and the cap.
\begin{proposition}\label{prop_bounds}
	For every fixed $k\geq3$, as $n\ra\infty$,
	\begin{align*}
		1.5-o_n(1)&\leq{\Erw\tau_n^+\over n}\leq2^k\ln2+o_n(1),\\
		1.514-o_n(1)&\leq{\Erw\tau_n^\pm\over n}\leq2^k\ln2+o_n(1).
	\end{align*}
\end{proposition}
\begin{proof}
	For either model, a fixed assignment violates a random clause with probability at least $2^{1-k}-O_k(1/n)$. Thus the first moment bound gives
	\begin{align*}
		\Pr(\Phi^\star_{n,\floor{2^k\ln2\cdot n}}\in\SAT)\leq\exp(-n\ln2+O_k(1)),\qquad\star\in\{+,\pm\}.
	\end{align*}
	Since $\tau_n^\star\leq2^k n$, this implies the upper bounds.
	For the lower bounds, extending each clause of a random NAE-3-SAT formula by $k-3$ uniformly chosen distinct variables outside the original clause, with independent uniform signs, preserves satisfiability and produces the NAE-$k$-SAT formula. The same construction with positive signs applies to hypergraph $2$-colourability.
	The bounds therefore follow from the results for NAE-3-SAT \cite{ACIM} and $3$-uniform hypergraph $2$-colourability \cite{AM2}.
\end{proof}

\section{Variance and concentration}

We begin with the NAE version of the clause replacement argument from \cite[Lemma 2.1]{AI}. We use the equivalence classes from \Sec~\ref{sec_classes}, rather than frozen variables of ordinary SAT.

\begin{proposition}\label{prop_short_clause}
	Fix $k\geq3$. There is $n_0\geq k+1$ such that the following holds for every $n\geq n_0$. Let $\Phi$ be any NAE-formula on $n-1$ variables, with arbitrary signs and clause lengths. For $\star\in\{+,\pm\}$, let $h^\star$ be a uniformly random $(k-1)$-clause and let $h_1^\star,\ldots,h_{g}^\star$ be independent uniformly random $k$-clauses, all on distinct variables within each clause, with signs of type $\star$, independent of one another. Then
	\begin{align}\label{eq_forcing_1}
		\Pr\bigl(\Phi\cup\{h_1^\star,\ldots,h_{g}^\star\}\notin\SAT\bigr)
		\geq \Pr\bigl(\Phi\cup\{h^\star\}\notin\SAT\bigr)-\eps,
	\end{align}
	where $\eps=n^{-(k-2)/(k-1)}$, $g=\ceil{k n^{1/(k-1)}}$.
\end{proposition}
\begin{proof}
	The assertion is immediate if $\Phi$ is unsatisfiable. Otherwise put $q=q_{k-1}^\star(\Phi)$. It suffices to consider $q>\eps$.
	By \Lem~\ref{lem_consecutive}, $q_k^\star(\Phi)\geq q^{(k-1)/(k-2)}-(k-2)q/(n-k)$.
	Since $q>n^{-(k-2)/(k-1)}$, for all sufficiently large $n$, uniformly in $\Phi$,
	\[
		\frac{g q_k^\star(\Phi)}q
		\geq g\left(n^{-1/(k-1)}-\frac{k-2}{n-k}\right)\geq2.
	\]
	Independence of the new clauses implies
	\[
		\Pr\bigl(\Phi\cup\{h_1^\star,\ldots,h_{g}^\star\}\notin\SAT\bigr)
		\geq1-(1-q_k^\star(\Phi))^{g}
		\geq1-e^{-g q_k^\star(\Phi)}
		\geq1-e^{-2q}.
	\]
	If $q\leq1/2$, then $1-e^{-2q}\geq q$, which proves the claim. If $q>1/2$, \Lem~\ref{lem_consecutive} bounds $q_k^\star(\Phi)$ away from zero uniformly in $\Phi$. Thus $e^{-g q_k^\star(\Phi)}\leq\eps$ for sufficiently large $n$, which gives a probability at least $1-\eps\geq q-\eps$.
\end{proof}

Using \Prop~\ref{prop_short_clause} we control the total expected number of forced variables along the clause process.
\begin{proposition}\label{prop_b_sum}
	Fix $k\geq3$ and $\lambda>0$. For every integer $0\leq M\leq\lambda n$ and either $\star\in\{+,\pm\}$,
	\begin{align*}
		\sum_{j=0}^M\Erw{b(\Phi^\star_{n,j})\over n}=O_{k,\lambda}(n^{1/(k-1)}).
	\end{align*}
\end{proposition}
\begin{proof}
	By symmetry, the $j$-th summand is the probability that $\Phi^\star_{n,j}$ is satisfiable and $x_n$ is forced. The number $\zeta$ of clauses containing $x_n$ has distribution $\Bin(j,k/n)$. Delete these clauses to obtain $\Phi'$, and erase the literal on $x_n$ from each deleted clause to obtain $h_1,\ldots,h_\zeta$.
	Conditional on the incidence positions, $\Phi'$ consists of $j-\zeta$ independent uniform $k$-clauses on the remaining $n-1$ variables. The $h_t$ are independent uniform $(k-1)$-clauses on these variables, independent of $\Phi'$, with signs of type $\star$.
	If $x_n$ is forced, then
	\begin{align*}
		\Phi'\in\SAT,\qquad\Phi'\cup\{h_1,\ldots,h_\zeta\}\notin\SAT.
	\end{align*}

	Put $P(m)=\Pr(\Phi^\star_{n-1,m}\in\SAT)$ and use $g$ and $\eps$ from \Prop~\ref{prop_short_clause}. Conditional on $\zeta=d$ and the incidence positions, replace each $h_t$ by $g$ fresh $k$-clauses.
	At each replacement, condition on all other clauses and apply \Prop~\ref{prop_short_clause}. The total error is at most $d\eps$, and subtracting the probability that $\Phi'$ is itself unsatisfiable gives
	\begin{align*}
		\Pr(x_n\text{ is forced}\mid\text{incidence positions})
		\leq P(j-d)-P(j-d+dg)+d\eps.
	\end{align*}
	For $0\leq d\leq j\leq M$, we have $\Pr(\zeta=d)\leq(kM/n)^d/d!$, and
	\begin{align*}
		\sum_{j=d}^M\bc{P(j-d)-P(j-d+dg)}\leq dg,
	\end{align*}
	where the sum telescopes and $P$ is nonincreasing. Averaging the error with $\Erw\zeta=kj/n$, we obtain
	\begin{align*}
		\sum_{j=0}^M\Erw{b(\Phi^\star_{n,j})\over n}
		&\leq g\sum_{d\geq0}{d(kM/n)^d\over d!}+\eps\sum_{j=0}^M{kj\over n}\\
		&={kM\over n}e^{kM/n}g+{kM(M+1)\over2n}\eps
		=O_{k,\lambda}(n^{1/(k-1)}).
	\end{align*}
\end{proof}

Fix an integer $M\geq1$ and $\star\in\{+,\pm\}$. Write $C_1^\star,\ldots,C_M^\star$ for the independent clauses of $\Phi^\star_{n,M}$ and set
\begin{align*}
	\tau^\star_n(M)=\max\cbc{0\leq m\leq M:~\Phi^\star_{n,m}\in\SAT}.
\end{align*}
For $1\leq i\leq M$, let $\Phi^\star_{m,-i}$ be the formula $\Phi^{\star}_{n,m}$ with the clause at position $i$ removed if $i\leq m$.
Set
\begin{align*}
	\tau^\star_{n,-i}(M)=\max\cbc{0\leq m\leq M:~\Phi^\star_{m,-i}\in\SAT}.
\end{align*}
In particular, $\tau^\star_{n,-i}(M)\geq\tau^\star_n(M)$, and $\tau^\star_{n,-i}(M)$ is independent of $C_i^\star$.

\begin{proposition}\label{prop_var_deletion}
	For $\star\in\{+,\pm\}$,
	\begin{align*}
		\Var(\tau^\star_n(M))\leq\sum_{i=1}^M\Erw(\tau^\star_{n,-i}(M)-\tau^\star_n(M))^2.
	\end{align*}
\end{proposition}
\begin{proof}
	Apply \cite[Eq. (3.6)]{BBLM} with $Z=\tau^\star_n(M)$ and $Z_i=\tau^\star_{n,-i}(M)$.
\end{proof}

\begin{proposition}\label{prop_var_summand}
	For every $1\leq i\leq M$ and either $\star\in\{+,\pm\}$,
	\begin{align*}
		\Erw(\tau^\star_{n,-i}(M)-\tau^\star_n(M))^2\leq2\sum_{m=i}^M\Erw\min\cbc{M q_k^\star(\Phi^\star_{m,-i}),1}.
	\end{align*}
\end{proposition}
\begin{proof}
	The proof of \cite[Lemma 3.2]{AI} applies with killing probability $q_k^\star$, as it only uses independence of the clauses and monotonicity of satisfiability under clause addition.
\end{proof}

\begin{lemma}\label{lem_var}
	For fixed $k\geq3$, $\lambda>0$, integers $1\leq M\leq\lambda n$ and $\star \in \{+, \pm\}$,
	\begin{align*}
		\Var(\tau^\star_n(M))=O_{k,\lambda}\bc{n^{1+1/k+1/(k-1)}}.
	\end{align*}
\end{lemma}
\begin{proof}
	By \eqref{eq_forced_bound},
	\begin{align*}
		\min\cbc{M q_k^\star(\Phi^\star_{m,-i}),1}\leq M^{1/k}{b(\Phi^\star_{m,-i})\over n}.
	\end{align*}
	For $m\geq i$, $\Phi^\star_{m,-i}$ has the same distribution as $\Phi^\star_{n,m-1}$. Thus \Prop s~\ref{prop_b_sum}--\ref{prop_var_summand} yield
	\begin{align*}
		\Var(\tau^\star_n(M))
		&\leq2M^{1/k}\sum_{i=1}^M\sum_{m=i}^M\Erw{b(\Phi^\star_{n,m-1})\over n}\\
		&\leq2M^{1+1/k}\sum_{j=0}^{M-1}\Erw{b(\Phi^\star_{n,j})\over n}
		=O_{k,\lambda}\bc{n^{1+1/k+1/(k-1)}}.
	\end{align*}
\end{proof}

\begin{corollary}\label{cor_concentration}
	Under the assumptions of \Lem~\ref{lem_var}, for either $\star\in\{+,\pm\}$,
	\begin{align*}
		\Pr\bc{|\tau^\star_n(M)-\Erw\tau^\star_n(M)|\geq n^{19/20}}=O_{k,\lambda}(n^{-1/15}).
	\end{align*}
\end{corollary}
\begin{proof}
	Apply Chebyshev's inequality and use $1+1/k+1/(k-1)\leq11/6$ for $k\geq3$.
\end{proof}

\section{From concentration to convergence}

We now use the concentration bound to prove that the normalized expectations converge.
For $\lambda>0$ and $\star\in\{+,\pm\}$, put
\begin{align*}
	M_{n,\lambda}&=\floor{\lambda n},&
	\tau^\star_{n,\lambda}&=\min\cbc{\xi_n^\star-1,M_{n,\lambda}},&
	\mu_{n,\lambda}^\star&={\Erw\tau_{n,\lambda}^\star\over n}.
\end{align*}
Thus $\tau_{n,\lambda}^\star=\tau^\star_n(M_{n,\lambda})$, and $\tau_{n,2^k}^\star=\tau_n^\star$.

\begin{theorem}\label{thm_bounded}
	Fix an integer $k\geq3$, $\lambda>0$ and $\star\in\{+,\pm\}$.
	Suppose there exist
	\begin{align*}
		0<a\leq b<\lambda,\qquad 0<\delta<1/2,\qquad\gamma>0,\qquad0\leq K<+\infty
	\end{align*}
	such that for all sufficiently large $n$,
	\begin{align}
		a\leq\mu^\star_{n,\lambda}\leq b,\label{eq_concentration_1}\\
		\Pr\bc{\abs{\tau^\star_{n,\lambda}-n\mu^\star_{n,\lambda}}\geq n^{1-\delta}}\leq Kn^{-\gamma}.\label{eq_concentration_2}
	\end{align}
	Then $\mu^\star_{n,\lambda}$ converges to $\alpha_k^\star\in[a,b]$, and
	\begin{align*}
		\lim_{n\ra\infty}\Pr(\tau^\star_{n,\lambda}\geq\floor{cn})
		=\begin{cases}
			1, & 0\leq c<\alpha_k^\star;\\
			0, & c>\alpha_k^\star.
		\end{cases}
	\end{align*}
\end{theorem}

\subsection{Poissonian model}

Fix $k\geq3$. For $n\geq1$ and $c\geq0$, let $\Psi^+_{n,cn}$ be the NAE-formula with $\Po(cn)$ independent positive clauses, each formed by sampling $k$ variables independently and uniformly from the $n$ variables. Thus repetitions within clauses are allowed.
Define $\Psi^\pm_{n,cn}$ in the same way, with independent uniform signs for each literal.
For $\star\in\{+,\pm\}$, write $s_n^\star(c)=\Pr(\Psi^\star_{n,cn}\in\SAT)$.

\begin{lemma}\label{lem_interpolation}
	For any $n_1,n_2\in\NN$, $n=n_1+n_2$, $c\geq0$ and $\star\in\{+,\pm\}$,
	\begin{align*}
		s^\star_n(c)\geq s^\star_{n_1}(c)\cdot s^\star_{n_2}(c).
	\end{align*}
\end{lemma}
\begin{proof}
	We use a Poisson version of the combinatorial interpolation of \cite[Proposition 2]{BGT}, as in \cite[Lemma 5.2]{AI}, with the killing probabilities from \Lem~\ref{lem_clause_classes}. Fix $\star\in\{+,\pm\}$ and partition the variable set $V$ into $V_1,V_2$ of sizes $n_1,n_2$. At time $t\in[0,1]$, sample independent numbers of global and local clauses
	\[
		\xi \sim \Po(cnt),\qquad \xi_1 \sim \Po(cn_1(1-t)),\qquad \xi_2 \sim \Po(cn_2(1-t)).
	\]
	Each clause has $k$ independently chosen variables from $V$, $V_1$ or $V_2$ respectively, and signs of type $\star$. Let $\Psi_t^\star$ be their union and $f^\star(t)=\Pr(\Psi_t^\star\in\SAT)$. Then
	\[
		f^\star(0)=s_{n_1}^\star(c)s_{n_2}^\star(c),\qquad
		f^\star(1)=s_n^\star(c).
	\]
	For each satisfiable realization of $\Psi_t^\star$, take the classes and parts from \Sec~\ref{sec_classes}. A clause confined to $V_h$ therefore uses the intersections $\cO_j\cap V_h$ (or $\cO_{j,a}\cap V_h$) for $h = 1,2$.
	The Poisson add-one identity says that differentiation with respect to a clause family's mean gives the expected change in the satisfiability indicator upon adding one clause, and this change is minus its killing probability. Therefore, \Lem~\ref{lem_clause_classes} gives
	\begin{align*}
		{\partial f^\pm(t) \over \partial t}
		&=cn\,2^{1-k}\Erw\brk{
		\vecone_{\{\Psi_t^\pm\in\SAT\}}\sum_j
		\bc{\sum_{h=1}^2{n_h\over n}\bc{|\cO_j\cap V_h|\over n_h}^k
		-\bc{|\cO_j|\over n}^k}}\geq0.
	\end{align*}
	Indeed, due to convexity,
	\begin{align*}
		\bc{|\cO_j|\over n}^k
		=\bc{\sum_{h=1}^2{n_h\over n}{|\cO_j\cap V_h|\over n_h}}^k
		\leq\sum_{h=1}^2{n_h\over n}\bc{|\cO_j\cap V_h|\over n_h}^k.
	\end{align*}
	For positive clauses, \Lem~\ref{lem_clause_classes} gives
	\begin{align*}
		{\partial f^+(t) \over \partial t}
		&=cn\Erw\brk{
		\vecone_{\{\Psi_t^+\in\SAT\}}\sum_{j=1}^J\sum_{a=0}^1
		\bc{\sum_{h=1}^2{n_h\over n}\bc{|\cO_{j,a}\cap V_h|\over n_h}^k
		-\bc{|\cO_{j,a}|\over n}^k}}\geq0.
	\end{align*}
	Integrating over $t$ proves both inequalities.
\end{proof}

\subsection{From the Poisson model to proper clauses}
We next transfer the concentration bounds to the Poisson model.
Fix $\star\in\{+,\pm\}$ and parameters $k,\lambda,a,b,\delta,\gamma,K$ from the hypothesis of \Thm~\ref{thm_bounded}.

\begin{lemma}\label{lem_transfer}
	Fix $k\geq3$ and $D\geq0$. For $n\geq k$ and $0\leq c\leq D$, put $\theta_n=(n)_k/n^k$ and let $\zeta \sim\Po(cn\theta_n)$. For either $\star\in\{+,\pm\}$,
	\begin{align}\label{eq_pois_01}
		e^{-D\binom{k}{2}}\Pr(\Phi^\star_{n,\zeta}\in\SAT)
		\leq s_n^\star(c)\leq\Pr(\Phi^\star_{n,\zeta}\in\SAT).
	\end{align}
	Moreover, $cn-D\binom{k}{2}\leq\Erw \zeta\leq cn$ and $\Var \zeta\leq Dn$. All bounds are uniform over $c\in[0,D]$.
\end{lemma}
\begin{proof}
	Apply the marked Poisson thinning argument of \cite[Lemma 5.3]{AI}. The numbers of clauses with distinct variables and with a repetition are independent with laws
	\[
		\zeta\sim\Po(cn\theta_n),\qquad \eta\sim\Po(cn(1-\theta_n)).
	\]
	Conditional on $\zeta$ and $\eta$, the two clause families are independent, and any clause with distinct variables is a uniformly random set of size $k$ with signs of type $\star$.
	Thus the proper family $\Psi'^{\star}$ has the law $\Phi^\star_{n,\zeta}$ in either model. With $\Psi^\star$ denoting the full formula,
	\begin{align}\label{eq_pois_1}
		\{\Psi'^\star\in\SAT\}\cap\{\eta=0\}
		\subseteq\{\Psi^\star\in\SAT\}\subseteq\{\Psi'^\star\in\SAT\}.
	\end{align}
	Independence and $\Pr(\eta=0)=e^{-cn(1-\theta_n)}\geq e^{-D\binom{k}{2}}$ give the two inequalities. The moment bounds follow from $\Erw \zeta=\Var \zeta=cn\theta_n$.
\end{proof}

\begin{lemma}\label{lem_impl}
	Assume \eqref{eq_concentration_1}--\eqref{eq_concentration_2} and fix $D>b$. For all sufficiently large $n$, uniformly over $0\leq c\leq D$,
	\begin{align}
		c\leq \mu^\star_{n,\lambda}-2n^{-\delta}
		&\quad\Longrightarrow\quad s_n^\star(c)\geq{1\over2}e^{-D\binom{k}{2}},\label{eq_pois_0_1}\\
		c\geq \mu^\star_{n,\lambda}+2n^{-\delta}
		&\quad\Longrightarrow\quad s_n^\star(c)=O\bc{n^{-\gamma}+n^{-(1-2\delta)}}.\label{eq_pois_0_2}
	\end{align}
\end{lemma}
\begin{proof}
	Take $\zeta$ from \Lem~\ref{lem_transfer} independently of the randomness inside clauses. By Chebyshev's inequality and \eqref{eq_concentration_2},
	\begin{align*}
		|\zeta-cn|\leq{1\over2}n^{1-\delta},\qquad
		|\tau^\star_{n,\lambda}-n\mu^\star_{n,\lambda}|<n^{1-\delta}
	\end{align*}
	except on an event of probability $O(n^{-\gamma}+n^{-(1-2\delta)})$ uniformly for $c\in[0,D]$.
	If $c\leq\mu^\star_{n,\lambda}-2n^{-\delta}$, these inequalities yield
	\begin{align*}
		\zeta\leq cn+{1\over2}n^{1-\delta}
		\leq n\mu^\star_{n,\lambda}-{3\over2}n^{1-\delta}
		<\tau^\star_{n,\lambda}\leq\xi_n^\star-1.
	\end{align*}
	Thus $\Phi^\star_{n,\zeta}$ is satisfiable outside the exceptional event, and \eqref{eq_pois_01} implies \eqref{eq_pois_0_1} for $n$ large enough.

	Conversely, if $c\geq\mu^\star_{n,\lambda}+2n^{-\delta}$, then
	\begin{align*}
		\zeta\geq cn-{1\over2}n^{1-\delta}
		\geq n\mu^\star_{n,\lambda}+{3\over2}n^{1-\delta}
		>\tau^\star_{n,\lambda}
	\end{align*}
	with probability $1 - O(n^{-\gamma} + n^{-(1-2\delta)})$.
	Outside the exceptional event, $\tau^\star_{n,\lambda}<bn+n^{1-\delta}<M_{n,\lambda}$ for large $n$, hence $\tau^\star_{n,\lambda}=\xi_n^\star-1$, so $\Phi^\star_{n,\zeta}$ is unsatisfiable. The upper bound in \eqref{eq_pois_01} now yields \eqref{eq_pois_0_2}.
\end{proof}

Fix $D>b$. The next lemma compares the expected threshold locations at different sizes.
\begin{lemma}\label{lem_super_add}
	For all sufficiently large $r,s$,
	\begin{align}
		\mu^\star_{r+s,\lambda}&\geq\min\cbc{\mu^\star_{r,\lambda},\mu^\star_{s,\lambda}}-4\min\{r,s\}^{-\delta},\label{eq_main_1}\\
		\mu^\star_{s+1,\lambda}&\geq \mu^\star_{s,\lambda}-4s^{-\delta}.\label{eq_main_2}
	\end{align}
\end{lemma}
\begin{proof}
	Choose $c=\min\cbc{\mu^\star_{r,\lambda},\mu^\star_{s,\lambda}}-2\min\{r,s\}^{-\delta}$. For large $r,s$ it lies in $[0,b]$. By \eqref{eq_pois_0_1} and \Lem~\ref{lem_interpolation},
	\begin{align*}
		s^\star_{r+s}(c)\geq s^\star_r(c)s^\star_s(c)\geq{1\over4}e^{-2D\binom{k}{2}}.
	\end{align*}
	Eq. \eqref{eq_pois_0_2} therefore implies $c<  \mu^\star_{r+s,\lambda}+2(r+s)^{-\delta}$ for large $r+s$, proving \eqref{eq_main_1}.

	For \eqref{eq_main_2}, take $c=\mu^\star_{s,\lambda}-2s^{-\delta}$. The one-variable formula is empty with probability $e^{-c}$, so $s^\star_1(c)\geq e^{-c}\geq e^{-D}$. Thus \Lem~\ref{lem_interpolation} and \eqref{eq_pois_0_1} imply that $s^\star_{s+1}(c)\geq s_s^\star(c)s_1^\star(c)\geq{1\over2}e^{-D(1+\binom{k}{2})}$, which is bounded away from zero. Applying \eqref{eq_pois_0_2} once more yields $\mu^\star_{s+1,\lambda}>c-2(s+1)^{-\delta}\geq \mu^\star_{s,\lambda}-4s^{-\delta}$.
\end{proof}

\subsection{Convergence}

\begin{lemma}\label{lem_tree}
	Let $(u_n)_{n \in \NN}$ be a bounded real-valued sequence. Suppose there exist $0 \leq A < +\infty$ and $\beta > 0$, such that for large enough $r, s$ with $1/3 \leq r/s \leq 3$,
	\begin{align*}
		u_{r + s} &\geq \min\cbc{u_r, u_s} - A \min\cbc{r,s}^{-\beta},
		\\
		u_{s + 1}& \geq u_s - A s^{-\beta}.
	\end{align*}
	Then $u_n$ converges.
\end{lemma}
\begin{proof}
	This is \cite[Lemma 5.4]{AI}.
\end{proof}

\subsection{Proof of \Thm~\ref{thm_bounded}}
Apply \Lem~\ref{lem_tree} to the sequence $\mu^\star_{n,\lambda}$ from \Lem~\ref{lem_super_add}, with $A=4$ and $\beta=\delta$. This gives $\mu^\star_{n,\lambda}\ra\alpha_k^\star\in[a,b]$.
For fixed $0\leq c<\alpha_k^\star$, the gap $\mu^\star_{n,\lambda}-c$ is bounded away from zero for large $n$, so \eqref{eq_concentration_2} implies $\tau^\star_{n,\lambda}\geq\floor{cn}$ \whp.
For fixed $c>\alpha_k^\star$, \eqref{eq_concentration_2} instead gives $\Pr(\tau^\star_{n,\lambda}\geq\floor{cn})\ra0$, completing the proof.

\subsection{Proof of \Thm~\ref{thm_main}}
Take $\lambda=2^k$. By \Prop~\ref{prop_bounds}, condition \eqref{eq_concentration_1} holds with $a=1$ and $b = 2^k - 1$. \Cor~\ref{cor_concentration} gives \eqref{eq_concentration_2} with $\delta=1/20$ and $\gamma=1/15$, thus \Thm~\ref{thm_bounded} can be applied. For $c<\lambda$ and large $n$, the events $\{\Phi^\star_{n,\floor{cn}}\in\SAT\}$ and $\{\tau^\star_{n,\lambda}\geq\floor{cn}\}$ coincide, and for $c\geq\lambda$ we use monotonicity of the property $\{\Phi^\star_{n,\floor{cn}} \in \SAT\}$.

\section{Threshold comparison}

In this section we prove \Thm~\ref{thm_comparison}.
For $c\geq0$ and $t\in[0,1]$, let $\Psi^{(t)}_{n,cn}$ be the random formula consisting of $\xi^+ \sim \Po(cnt)$ positive clauses and $\xi^{\pm} \sim \Po(cn(1-t))$ signed clauses, with $\xi^+$ and $\xi^\pm$ being independent. Each clause samples its $k$ variables independently and uniformly, with repetitions allowed, and signs in all $\xi^{\pm}$ clauses are independent and uniformly random.

For each satisfiable $\Phi$, use the classes $\cO_j$ and parts $\cO_{j,a}$ from \Sec~\ref{sec_classes}. By \Lem~\ref{lem_clause_classes},
\[
	\widehat q_k^+(\Phi)=\sum_j\sum_{a=0}^1\left(\frac{|\cO_{j,a}|}{n}\right)^k,
	\qquad
	\widehat q_k^\pm(\Phi)=2^{1-k}\sum_j\left(\frac{|\cO_j|}{n}\right)^k,
\]
and convexity implies $\widehat q_k^+(\Phi)\geq\widehat q_k^\pm(\Phi)$ for $k \geq 3$.

Writing $f_n(t,c)=\Pr(\Psi^{(t)}_{n,cn}\in\SAT)$, the same Poisson add-one identity as in \Lem~\ref{lem_interpolation} yields
\[
	\frac{\partial f_n(t,c)}{\partial t}
	=cn\Erw\left[\vecone_{\{\Psi^{(t)}_{n,cn}\in\SAT\}}
	\left(\widehat q_k^\pm(\Psi^{(t)}_{n,cn})-\widehat q_k^+(\Psi^{(t)}_{n,cn})\right)\right]\leq0.
\]
Hence $s_n^+(c)\leq s_n^\pm(c)$ for every $n$ and $c\geq0$.
Suppose that for some $k$, $\alpha_k^{\mathrm{NAE}}<c<\alpha_k^{\mathrm{2-col}}$. \Thm~\ref{thm_main}, Poisson concentration and monotonicity imply that the distinct-variable formula with $\zeta\sim\Po(cn\theta_n)$ positive clauses is satisfiable with probability tending to one, while the corresponding signed formula is satisfiable with probability tending to zero. Here $\theta_n=(n)_k/n^k$ and $\zeta/n$ converges in probability to $c$.
Furthermore, the coupling from \Lem~\ref{lem_transfer} implies that $s_n^+(c)$ is bounded away from zero, and that $s_n^\pm(c)\to0$, which leads to a contradiction.
Thus $\alpha_k^{\mathrm{2-col}}\leq\alpha_k^{\mathrm{NAE}}$.


\begin{thebibliography}{99}
	\bibitem{AbMo} E.~Abbe, A.~Montanari: On the Concentration of the Number of Solutions of Random Satisfiability Formulas. Random Structures \& Algorithms 45 (2014) 362--382.
	\bibitem{ACIM} D.~Achlioptas, A.~Chtcherba, G.~Istrate, C.~Moore: The Phase Transition in 1-in-$k$ SAT and NAE 3-SAT. Proceedings of the 12th Annual ACM-SIAM Symposium on Discrete Algorithms (2001) 721--722.
	\bibitem{AF} D.~Achlioptas, E.~Friedgut: A Sharp Threshold for $k$-Colorability. Random Structures \& Algorithms 14 (1999) 63--70.
	\bibitem{AM} D.~Achlioptas, C.~Moore: On the 2-Colorability of Random Hypergraphs. Randomization and Approximation Techniques in Computer Science (2002) 78--90.
	\bibitem{AM2} D.~Achlioptas, C.~Moore: Random $k$-SAT: Two Moments Suffice to Cross a Sharp Threshold. SIAM Journal on Computing 36 (2006) 740--762.
	\bibitem{AN} D.~Achlioptas, A.~Naor: The Two Possible Values of the Chromatic Number of a Random Graph. Annals of Mathematics 162 (2005) 1335--1351.
	\bibitem{BGT} M.~Bayati, D.~Gamarnik, P.~Tetali: Combinatorial Approach to the Interpolation Method and Scaling Limits in Sparse Random Graphs. Annals of Probability 41 (2013) 4080--4115.
	\bibitem{BBLM} S.~Boucheron, O.~Bousquet, G.~Lugosi, P.~Massart: Moment Inequalities for Functions of Independent Random Variables. Annals of Probability 33 (2005) 514--560.
	\bibitem{Car2} G.~Carenini: The Scaling Window of Random $k$-SAT. arXiv:2609.26222 (2026).
	\bibitem{Car} G.~Carenini: A Polynomial Scaling Window for Random $k$-SAT and a Proof of the Satisfiability Conjecture. Electronic Colloquium on Computational Complexity, Report No.~229 (2026).
	\bibitem{CR} V.~Chv\'atal, B.~Reed: Mick Gets Some (the Odds Are on His Side). Proceedings of the 33rd Annual Symposium on Foundations of Computer Science (1992) 620--627.
	\bibitem{CIKZ} A.~Coja-Oghlan, M.~Isaev, M.~Kang, P.~Zakharov: The Chromatic Number of a Random Graph with $o(n^{3/2})$ Edges. arXiv:2610.33678 (2026).
	\bibitem{CP} A.~Coja-Oghlan, K.~Panagiotou: Catching the $k$-NAESAT Threshold. Proceedings of the 44th Annual ACM Symposium on Theory of Computing (2012) 899--908.
	\bibitem{CZ} A.~Coja-Oghlan, L.~Zdeborov\'a: The Condensation Transition in Random Hypergraph 2-Coloring. Proceedings of the 23rd Annual ACM-SIAM Symposium on Discrete Algorithms (2012) 241--250.
	\bibitem{DSS} J.~Ding, A.~Sly, N.~Sun: Proof of the Satisfiability Conjecture for Large $k$. Annals of Mathematics 196 (2022) 1--388.
	\bibitem{ER} P.~\Erdos, A.~\Renyi: On the Evolution of Random Graphs. Publications of the Mathematical Institute of the Hungarian Academy of Sciences 5 (1960) 17--61.
	\bibitem{Fri} E.~Friedgut, with an appendix by J.~Bourgain: Sharp Thresholds of Graph Properties, and the $k$-SAT Problem. Journal of the American Mathematical Society 12 (1999) 1017--1054.
	\bibitem{FriH} E.~Friedgut: Hunting for Sharp Thresholds. Random Structures \& Algorithms 26 (2005) 37--51.
	\bibitem{Heckel} A.~Heckel: The Chromatic Number of Dense Random Graphs. Random Structures \& Algorithms 53 (2018) 140--182.
	\bibitem{KMRTZ} F.~Krzakala, A.~Montanari, F.~Ricci-Tersenghi, G.~Semerjian, L.~Zdeborov\'a: Gibbs States and the Set of Solutions of Random Constraint Satisfaction Problems. Proceedings of the National Academy of Sciences 104 (2007) 10318--10323.
	\bibitem{AI} OpenAI: A Limiting Satisfiability Threshold for Every Fixed Clause Size. preprint (2026).
\end{thebibliography}
\end{document}